\documentclass[a4paper,12pt]{amsart}

\usepackage{amsmath,amssymb,xcolor,verbatim,microtype}
\usepackage[top=3.5cm,bottom=3.5cm,left=2.5cm,right=2.5cm]{geometry}
\usepackage[linktoc=page,colorlinks,linkcolor={red!80!black},citecolor={red!80!black},urlcolor={blue!80!black}]{hyperref}\newcommand{\arxiv}[1]{\href{http://arxiv.org/pdf/#1}{arXiv:#1}}
\usepackage{tikz}
\usetikzlibrary{matrix,arrows}
\pgfrealjobname{homogeneous}

\usepackage{mathptmx}                                                              
\DeclareSymbolFont{cmletters}{OML}{cmm}{m}{it}                                     
\DeclareSymbolFont{cmlargesymbols}{OMX}{cmex}{m}{n}
\DeclareMathSymbol{\mycoprod}{\mathop}{cmlargesymbols}{"60}
\DeclareMathSymbol{\mydelta}{\mathord}{cmletters}{"0E}     \let\delta\mydelta
\DeclareMathSymbol{\myepsilon}{\mathord}{cmletters}{"0F}   \let\epsilon\myepsilon
\DeclareMathSymbol{\myeta}{\mathord}{cmletters}{"11}       \let\eta\myeta
\DeclareMathSymbol{\myiota}{\mathord}{cmletters}{"13}      \let\iota\myiota
\DeclareMathSymbol{\mylambda}{\mathord}{cmletters}{"15}    \let\lambda\mylambda
\DeclareMathSymbol{\mymu}{\mathord}{cmletters}{"16}        \let\mu\mymu
\DeclareMathSymbol{\mypi}{\mathord}{cmletters}{"19}        \let\pi\mypi

\theoremstyle{plain}
\newtheorem*{problem*}{Problem}
\newtheorem*{thm*}{Theorem}
\newtheorem*{cor*}{Corollary}

\newcommand\Hom{{\mathrm {Hom}}}

\newcommand\C{{\mathbb C}}
\renewcommand\P{{\mathbb P}}
\newcommand\cB{{\mathcal B}}

\newcommand\ud{{\underline{d}}}
\newcommand\ue{{\underline{e}}}
\newcommand\udim{{\underline{\dim}\, }}

\DeclareMathOperator{\Gr}{Gr}

\DeclareMathOperator{\Spec}{Spec}
\DeclareMathOperator{\Rep}{Rep}

\newenvironment{psmallmatrix}{\left[\begin{smallmatrix}}{\end{smallmatrix}\right]}

\begin{document}

\title{Pl\"ucker relations for quiver Grassmannians}

\author{Oliver Lorscheid}
\address{Instituto Nacional de Matem\'atica Pura e Aplicada, Rio de Janeiro, Brazil}
\email{\href{mailto:oliver@impa.br}{oliver@impa.br}}

\author{Thorsten Weist}
\address{Bergische Universit\"at Wuppertal, Gau\ss str.\ 20, 42097 Wuppertal, Germany}
\email{\href{mailto:weist@uni-wuppertal.de}{weist@uni-wuppertal.de}}

\begin{abstract}
 In this text, we exhibit the \emph{quiver Pl\"ucker relations} for a quiver Grassmannian and show that they describe the quiver Grassmannian as a closed subscheme of a product of usual Grassmannians.
\end{abstract}

\maketitle

\subsection*{Introduction}
Let $Q$ be a quiver, $M$ be a finite dimensional complex representation of $Q$ and $\ue=(e_p)_{p\in Q_0}$ a dimension vector of $Q$. As a set, the quiver Grassmannian $\Gr_\ue(M)$ can be defined as the collection of all subrepresentations $N$ of $M$ with dimension vector $\udim N=\ue$. It gains the structure of a complex projective variety in terms of the closed embedding
\[
 \iota: \ \Gr_\ue(M) \ \longrightarrow \ \prod_{p\in Q_0} \Gr(e_p,M_p) 
\]
into the product of usual Grassmannians $\Gr(e_p,M_p)$ which sends a subrepresentation $N$ of $M$ to the collection of $e_p$-dimensional subvector spaces $N_p$ of $M_p$. 

The choice of ordered bases $\cB_p$ of $M_p$ for each vertex $p$ of $Q$ identifies $M_p$ with $\C^{d_p}$ and yields Pl\"ucker coordinates $[:\Delta_I:]$ for the complex points of $\Gr(e_p,M_p)=\Gr(e_p,d_p)$ where $I$ varies through all cardinality $e_p$-subsets of $\cB_p$. The following question suggests itself.

\begin{problem*}
 Can we exhibit explicit equations for the image of $\Gr_\ue(M)$ in $\prod \Gr(e_p,d_p)$?
\end{problem*}

This question has not been answered yet, possibly for the following reasons. First of all, quiver Grassmannians fail to be homogeneous spaces in general, which means that the shape of the defining equations depend on the choice of the ordered basis $\cB$. Secondly, $\Gr_\ue(M)$ carries a schematic structure coming from its description as a fibre of the universal quiver Grassmannian of $Q$, and this structure as a scheme is not visible from the viewpoint of a pointwise embedding into $\prod\Gr(e_p,d_p)$. Finally, we note that the situation is as bad as possible: every projective scheme is isomorphic to a quiver Grassmannian; cf.\ \cite{Hille15} or \cite{Reineke12}.

\subsection*{Quiver Pl\"ucker relations}
Let $m_{v,j,i}$ be the matrix coefficients of the linear maps $M_v:M_p\to M_q$. For a subset $I\cup\{i\}$ of $\cB_p$, we define $\epsilon(i,I)=\#\{i'\in I|i'\leq i\}$. In this text, we show that the image of $\Gr_\ue(M)$ in $\prod\Gr(e_p,M_p)$ is carved out by the \emph{quiver Pl\"ucker relations}
\[
 E(v,I,J) \ = \ \sum_{i\in \cB_p-I,\, j\in J} (-1)^{\epsilon(i,I)+\epsilon(j,J)} \, m_{v,j,i} \, \Delta_{I\cup\{i\}} \, \Delta_{J-\{j\}} \ = \ 0
\]
where $v:p\to q$ varies through the arrows of $Q$, $I$ varies through the $(e_p-1)$-subsets of $\cB_p$ and $J$ varies through the $(e_q+1)$-subsets of $\cB_q$.

In fact we prove the following stronger result: the quiver Pl\"ucker relations describe the structure of $\Gr_\ue(M)$ as a closed subscheme of $\prod \Gr(e_p,M_p)$. 

\subsection*{Remark on the changes to a previous version}
This paper is an improved version of what has been titled \emph{Homogeneous coordinates for quiver Grassmannians} before. We were able to remove the restrictive hypothesis on a dense Schubert cell from the main result of the previous version. As a result, we obtain the unconditional theorem in this text. Since both the formulas and the proof are different and various remarks are obsolete, we rewrote the text from scratch.

\subsection*{The quiver Grassmannian as a scheme}
We follow the exposition in \cite{cr}. Let $\ud=\udim M$ be the dimension vector of $M$. Define
\[
 \Rep(\ud) \ = \ \prod_{v:p\to q\text{ in }Q} \Hom(\C^{d_p},\C^{d_q}),
\]
which is the moduli space of $\ud$-dimensional representations of $Q$ with fixed basis. The complex points of $\Rep(\ud)$ correspond to tuples of matrices $A_v\in\Hom(\C^{d_p},\C^{d_q})$,
indexed by the arrows $v:p\to q$ of $Q$. Consider the trivial fibre bundle 
\[
 \pi: \ \Rep(\ud)\, \times\prod_{p\in Q_0}\Gr(e_p,d_p) \ \longrightarrow \ \Rep(\ud)
\]
whose complex points correspond to pairs of a tuple $(A_v)_{v\in Q_1}$ of matrices $A_v\in\Hom(\C^{d_p},\C^{d_q})$ and a tuple $(N_p)_{p\in Q_0}$ of $e_p$-dimensional subspaces $N_p$ of $\C^{d_p}$. 

The \emph{universal quiver Grassmannian} is the closed subscheme $\Gr_\ue^Q(\ud)$ of $\Rep(\ud)\times\prod\Gr(e_p,d_p)$ that is defined by the conditions $A_v(N_p) \subset N_q$ for every arrow $v:p\to q$ in $Q$. Let $\pi_Q$ be the restriction of $\pi$ to $\Gr_\ue^Q(\ud)$.

The choice of ordered bases $\cB_p$ of $M_p$ defines matrices $A_v=(m_{v,j,i})_{i\in\cB_p,j\in\cB_q}$ and thus a point $(A_v)_{v\in Q_1}$ of $\Rep(\ud)$. As a scheme, the quiver Grassmannian $\Gr_\ue(M)$ is defined as the fibre of $\pi_Q$ over $(A_v)_{v\in Q_1}$. Note that the embedding of $\Gr_\ue(M)$ into the fibre of $\pi$ over $(A_v)$ coincides with the embedding into $\prod\Gr(e_p,d_p)$ considered above.

\subsection*{Remark on the choice of bases}
A priori, the definition of $\Gr_\ue(M)$ as a complex scheme depends on the choice of the bases $\cB_p$, but it is easily seen that a change of bases induces a canonical isomorphism between the corresponding schemes. In so far, we can consider $\Gr_\ue(M)$ abstractly as a scheme. This will be, however, of no relevance for this text since we will work with a fixed choice of ordered bases $\cB_p$.

\subsection*{Theorem}\it
The image of the closed embedding 
\[
 \iota:\Gr_\ue(M)\to \prod\Gr(e_p,d_p)
\]
is the subscheme defined by the quiver Pl\"ucker relations
\[
 E(v,I,J) \ = \ \sum_{i\in \cB_p-I,\, j\in J} (-1)^{\epsilon(i,I)+\epsilon(j,J)} \, m_{v,j,i} \, \, \Delta_{I\cup\{i\}} \Delta_{J-\{j\}} \ = \ 0
\]
where $v:p\to q$ ranges through all arrows of $Q$, $I$ ranges through all $(e_p-1)$-subsets of $\cB_p$ and $J$ ranges through all $(e_q+1)$-subsets of $\cB_q$.
\rm

\subsection*{Proof}
The quiver Grassmannian is defined as a closed subscheme of $\Rep(\ud)\times\prod \Gr(e_p,d_p)$ by the relations
\[
 A_v \ = \ (m_{v,j,i})_{i\in \cB_p,j\in\cB_q} \qquad \text{and} \qquad A_v(N_p) \ \subset \ N_q
\]
for every arrow $v:p\to q$ in $Q$. The main effort of this proof is to show that the relations $A_v(N_p) \subset N_q$ are equivalent to the quiver Pl\"ucker relations.

To begin with, we remark that we can consider these relations for every arrow $v$ separately. This means that we can fix $v:p\to q$ for the following considerations. 

Consider the subspaces $N_p$ of $M_p$ and $N_q$ of $M_q$ with respective Pl\"ucker coordinates $[:\Delta_{I'}:]$ and $[:\Delta_{J'}:]$ where $I'$ varies through the $e_p$-subsets of $\cB_p$ and $J'$ varies through the $e_q$-subsets of $\cB_q$. We fix an $e_p$-subset $I_0$ of $\cB_p$ and an $e_q$-subset $J_0$ of $\cB_q$ such that $\Delta_{I_0}\neq0$ and $\Delta_{J_0}\neq0$.

Then $N_p$ is spanned by the vectors $n_{i_0}=(n_{i,i_0})_{i\in \cB_p}$ for $i_0\in I_0$ where
\[
 n_{i,i_0} \ = \ \begin{cases} \delta_{i,i_0} & \text{if }i\in I_0, \\
                               (-1)^{\epsilon(i,I)+\epsilon(i_0,I)}\, \frac{\Delta_{I\cup\{i\}}}{\Delta_{I_0}} & \text{if }i\notin I_0
                 \end{cases}
\]
with $I=I_0-\{i_0\}$. Note that $I$ is of cardinality $e_p-1$. 
Similarly, $N_q$ is spanned by the vectors $n_{j}=(n_{j_0,j})_{j_0\in \cB_q}$ for $j\in J_0$ where
\[
 n_{j_0,j} \ = \ \begin{cases} \delta_{j_0,j} & \text{if }j_0\in J_0, \\
                               (-1)^{\epsilon(j_0,J)+\epsilon(j,J)+1}\, \frac{\Delta_{J-\{j\}}}{\Delta_{J_0}} & \text{if }j_0\notin J_0
                 \end{cases}
\]
with $J=J_0\cup\{j_0\}$. Note that $J$ is of cardinality $e_q+1$. 

We conclude that $A_v(N_p)\subset N_q$ if and only if $A_v(n_{i_0})\in N_q$ for all $i_0\in I_0$, i.e.\ if there are $\lambda_{i_0,j_0}\in\C$ such that
\[
 A_v(n_{i_0}) \ = \ \sum_{j_0\in J_0} \lambda_{i_0,j_0}\, n_{j_0}.
\]
Since $n_{j_0,j}=\delta_{j_0,j}$ for $j\in J_0$, we have
\[
 \lambda_{i_0,j} \ = \sum_{j_0\in J_0} \lambda_{i_0,j_0}\, n_{j_0,j} \ = \ A_v(n_{i_0})_j \ = \ \sum_{i\in \cB_p} m_{v,j,i}\, n_{i,i_0}
\]
for $j\in J_0$. Thus $A_v(N_p)\subset N_q$ if and only if for all $i_0\in I_0$ and all $j_0\in \cB_q-J_0$,
\[
 \sum_{i\in\cB_p} m_{v,j_0,i}\, n_{i,i_0} \ = \ A_v(n_{i_0})_{j_0} \ = \ \sum_{j\in J_0} \lambda_{i_0,j}\, n_{j_0,j} \ = \ \sum_{j\in J_0,i\in \cB_p} m_{v,j,i}\, n_{i,i_0}\, n_{j_0,j}.
\]
Since $n_{i,i_0}=\delta_{i,i_0}$ for $i\in I_0$, this latter equation becomes
\[
 \sum_{j\in J_0,i\in \cB_p-I_0} m_{v,j,i}\, n_{i,i_0}\, n_{j_0,j} + \sum_{j\in J_0} m_{v,j,i_0}\, n_{j_0,j} - \sum_{i\in \cB_p-I_0} m_{v,j_0,i}\, n_{i,i_0} - m_{v,j_0,i_0} \ = \ 0.
\]
Replacing $n_{i,i_0}=(-1)^{\epsilon(i,I)+\epsilon(i_0,I)}\, \frac{\Delta_{I\cup\{i\}}}{\Delta_{I_0}}$ and $n_{j_0,j}=(-1)^{\epsilon(j_0,J)+\epsilon(j,J)+1}\, \frac{\Delta_{J-\{j\}}}{\Delta_{J_0}}$, this becomes
\begin{multline*}
 \sum_{j\in J_0,i\in \cB_p-I_0} (-1)^{\epsilon(i,I)+\epsilon(i_0,I)+\epsilon(j_0,J)+\epsilon(j,J)+1} \ m_{v,j,i} \ \frac{\Delta_{I\cup\{i\}}}{\Delta_{I_0}} \ \frac{\Delta_{J-\{j\}}}{\Delta_{J_0}} \\
 + \sum_{j\in J_0} (-1)^{\epsilon(j_0,J)+\epsilon(j,J)+1}  m_{v,j,i_0}  \frac{\Delta_{J-\{j\}}}{\Delta_{J_0}} \
 - \!\!\!\!\! \sum_{i\in \cB_p-I_0} (-1)^{\epsilon(i,I)+\epsilon(i_0,I)}  m_{v,j_0,i}  \frac{\Delta_{I\cup\{i\}}}{\Delta_{I_0}}  \
 - \ m_{v,j_0,i_0} \ = \ 0.
\end{multline*}
where $I=I_0-\{i_0\}$ and $J=J_0\cup\{j_0\}$. Multiplying this equation by $(-1)^{\epsilon(i_0,I)+\epsilon(j_0,J)+1}\Delta_{I_0}\Delta_{J_0}$ yields the quiver Pl\"ucker relation
\[
 E(v,I,J) \ = \ \sum_{i\in \cB_p-I,\, j\in J} (-1)^{\epsilon(i,I)+\epsilon(j,J)} \, m_{v,j,i} \, \, \Delta_{I\cup\{i\}} \Delta_{J-\{j\}} \ = \ 0.
\]

This shows that $A_v(N_p)\subset N_q$ if and only if for all $i_0\in I_0$ and all $j_0\in\cB_q-J_0$, we have $E(v,I,J)=0$ where $I=I_0-\{i_0\}$ and $J=J_0\cup\{j_0\}$. We conclude that the quiver Pl\"ucker relations define a subscheme of $\Gr_\ue(M)$.

We continue with the proof that this subscheme is indeed equal to $\Gr_\ue(M)$. Consider the quiver Pl\"ucker relations $E(v,I,J)=0$ for some choice of $v$, $I$ and $J$. If $\Delta_{I\cup\{i\}}=0$ for all $i\in\cB_p-I$ or $\Delta_{J-\{j\}}=0$ for all $j\in J$, then $E(v,I,J)=0$ is a trivial relation. Thus we can assume that there exist an $i\in \cB_p-I$ and a $j\in J$ such that $\Delta_{I_0}\neq0$ for $I_0=I\cup\{i_0\}$ and $\Delta_{J_0}\neq0$ for $J_0=J-\{j_0\}$. But for this choice of $I_0$ and $J_0$, our above considerations show that the quiver Pl\"ucker relation $E(v,I,J)=0$ appears among those equations characterizing the condition $A_v(N_p)\subset N_q$. This concludes the proof of the theorem. \qed

\subsection*{Remark on higher order and classical Pl\"ucker relations}
For every path $\pi$ from $p$ to $q$ along arrows $v_1,\dotsc,v_n$ in $Q$, we obtain the map $M_\pi=M_{v_n}\circ\dotsb\circ M_{v_1}:M_p\to M_q$. We call $n$ the \emph{length of $\pi$}. If the length of $\pi$ is $0$, i.e.\ if $\pi$ is the trivial path from $p$ to $q=p$, then we define $M_\pi$ as the identity map on $M_p=M_q$. Let $m_{\pi,j,i}$ denote the matrix coefficients of $M_\pi$ with respect to the bases $\cB_p$ of $M_p$ and $\cB_q$ of $M_q$.

As a consequence of the proof of the theorem, we see that the Pl\"ucker coordinates of the quiver Grassmannian $\Gr_\ue(M)$ satisfy the $n$-th order quiver Pl\"ucker relations
\[
 E(\pi,I,J) \ = \ \sum_{i\in \cB_p-I,\, j\in J} (-1)^{\epsilon(i,I)+\epsilon(j,J)} \, m_{\pi,j,i} \, \Delta_{I\cup\{i\}} \Delta_{J-\{j\}} \ = \ 0
\]
for every path $\pi$ from $p$ to $q$ of length $n$, every $(e_p-1)$-subset $I$ of $\cB_p$ and every $(e_q+1)$-subset $J$ of $\cB_q$. Note that the zeroth order quiver Pl\"ucker relations are nothing else than the classical Pl\"ucker relations
\[
 \sum_{i\in J-I} (-1)^{\epsilon(i,I)+\epsilon(i,J)} \, \Delta_{I\cup\{i\}} \Delta_{J-\{i\}} \ = \ 0.
\]

\subsection*{Remark on Schubert cells}
We can derive the equations from \cite[section 1.3]{LW15a} for a Schubert cell of $\Gr_\ue(M)$ from the quiver Pl\"ucker relations as follows: first we force the appropriate set of Pl\"ucker coordinates to be zero in the quiver Pl\"ucker relations $E(v,I,J)$, and then the equations for the Schubert cell result from dehomogenizing these equations.

\subsection*{Example 1 (Del Pezzo surface)}
The following quiver Grassmannian $\Gr_\ue(M)$ has already been investigated in Example 4.7 in \cite{L15}. The aim of our discussion is to illustrate how we can use the methods of this paper to find an explicit description of $\Gr_\ue(M)$. 

Let $Q$ be the Dynkin quiver of type $D_4$ in subspace orientation with arrows $a$, $b$ and $c$. Let $M$ be the representation that is given by the following coefficient quiver $\Gamma=\Gamma(M,\cB)$:
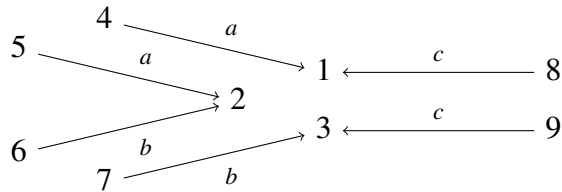
\[
  \beginpgfgraphicnamed{tikz/fig1}
  \begin{tikzpicture}[description/.style={fill=white,inner sep=0pt}]
   \matrix (m) [matrix of math nodes, row sep=-0.2em, column sep=1.5em, text height=1ex, text depth=0ex]
    {   & 4 &   &   &   &   &   &   &   \\
      5 &   &   &   &   &   &   &   &   \\     
        &   &   &   &   &   &   &   &   \\     
        &   &   &   & 1 &   &   &   & 8  \\     
        &   &   & 2 &   &   &   &   &   \\     
        &   &   &   & 3 &   &   &   & 9  \\     
        &   &   &   &   &   &   &   &   \\     
      6 &   &   &   &   &   &   &   &   \\     
        & 7 &   &   &   &   &   &   &   \\ };
    \path[->,font=\scriptsize]
    (m-1-2) edge node[auto] {$a$} (m-4-5)
    (m-2-1) edge node[auto] {$a$} (m-5-4)
    (m-4-9) edge node[auto,swap] {$c$} (m-4-5)
    (m-6-9) edge node[auto,swap] {$c$} (m-6-5)
    (m-8-1) edge node[auto,swap] {$b$} (m-5-4)
    (m-9-2) edge node[auto,swap] {$b$} (m-6-5);
  \end{tikzpicture}
  \endpgfgraphicnamed
\]
In other words, $M$ is given by the three matrices
\[
 A_a \ = \ \begin{psmallmatrix}
            1 & 0 \\
            0 & 1 \\
            0 & 0 \\
           \end{psmallmatrix}, 
     \qquad
 A_b \ = \ \begin{psmallmatrix}
            0 & 0 \\
            1 & 0 \\
            0 & 1 \\
           \end{psmallmatrix}, 
     \qquad
 A_c \ = \ \begin{psmallmatrix}
            1 & 0 \\
            0 & 0 \\
            0 & 1 \\
           \end{psmallmatrix}.
\]
Consider the dimension vector $\ue=(2,1,1,1)$ where the first coordinate refers to the central vertex of $Q$. The product Grassmannian $\prod\Gr(e_q,d_q)$ has Pl\"ucker coordinates
\[
 \bigr[ \; \Delta_{12}:\Delta_{13}:\Delta_{23} \; \bigr| \; \Delta_4:\Delta_5 \; \bigr| \; \Delta_6:\Delta_7 \; \bigr| \; \Delta_8:\Delta_9 \; \bigr]
\]
where we use the shorthand notations $\Delta_i=\Delta_{\{i\}}$ and $\Delta_{ij}=\Delta_{\{i,j\}}$. Since $\Gr(2,3)\simeq\P^2$ and $\Gr(1,2)=\P^1$, there are no classical Pl\"ucker relations among the coordinates. Thus the only relations between the coordinates are the quiver Pl\"ucker relations
\begin{align*}
 E(a,\emptyset,\{1,2,3\}) \ = \ \Delta_5 \; \Delta_{13} \ - \ \Delta_4 \; \Delta_{23} \ = \ 0, \\
 E(b,\emptyset,\{1,2,3\}) \ = \ \Delta_6 \; \Delta_{13} \ - \ \Delta_7 \; \Delta_{12} \ = \ 0, \\
 E(c,\emptyset,\{1,2,3\}) \ = \ \Delta_9 \; \Delta_{12} \ - \ \Delta_8 \; \Delta_{23} \ = \ 0. 
\end{align*}
It follows that the homogeneous coordinates $[\Delta_{12}:\Delta_{13}:\Delta_{23}]$ are determined by the other coordinates, and thus the embedding $\Gr_\ue(M)\to\prod\Gr(e_q,d_q)$ followed by the projection to the coordinates $\Delta_4,\dotsc,\Delta_9$ defines a closed embedding $\Gr_\ue(M)\to\P^1\times\P^1\times\P^1$.

For understanding the relation between the coordinates $\Delta_4,\dotsc,\Delta_9$, we consider the cases where one of $\Delta_{12}$, $\Delta_{13}$ and $\Delta_{23}$ is invertible. If $\Delta_{12}$ is invertible, then necessarily also $\Delta_6$ and $\Delta_8$ are invertible, and we obtain from the last two equations that
\[
 \Delta_{13} \ = \ \Delta_7 \; \Delta_{12} \; \Delta_6^{-1} \qquad \text{and} \qquad \Delta_{23} \ = \ \Delta_9 \; \Delta_{12} \; \Delta_8^{-1} 
\]
Substituting these terms in the first equation and multiplying with $\Delta_6\Delta_8\Delta_{12}^{-1}$ yields
\[
 \Delta_5 \; \Delta_7 \; \Delta_8 \ = \ \Delta_4 \; \Delta_6 \; \Delta_9. 
\]
Having $\Delta_{13}$ or $\Delta_{23}$ invertible results in the same equation. Thus $\Gr_\ue(M)$ is an irreducible surface of multidegree $(1,1,1)$ in $\P^1\times\P^1\times\P^1$. As argued in Example 4.7 of \cite{L15}, one finds that this surface is a del Pezzo surface of degree $6$.

\subsection*{Example 2 (Family with jumping Euler characteristic)}
The following example is taken from \cite{LW17}. Let $Q$ be a quiver of extended Dynkin quiver type $\widetilde A_2$ with arrows $a$, $b$ and $c$ as illustrated below. Let $\lambda\in\C$ be a parameter. Let $M_\lambda$ be the representation of $Q$ that is given by the following coefficient quiver $\Gamma_\lambda$ and the map $F:\Gamma_\lambda \to Q$. 
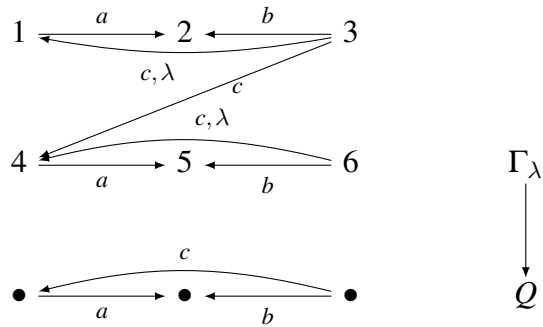
 \[
  \beginpgfgraphicnamed{tikz/fig2}
  \begin{tikzpicture}[>=latex]
   \matrix (m) [matrix of math nodes, row sep=3em, column sep=4em, text height=1ex, text depth=0ex]
    {     1  &    2    &    3                              \\   
          4  &    5    &    6    & \Gamma_\lambda   \\   
     \bullet & \bullet & \bullet & Q                \\   
    };
   \path[->,font=\scriptsize]
    (m-1-1) edge node[auto] {$a$} (m-1-2)
    (m-1-3) edge node[auto,swap] {$b$} (m-1-2)
    (m-1-3) edge[bend left=10] node[below,pos=0.6] {$c,\lambda$} (m-1-1)
    (m-1-3) edge node[right=0.3,pos=0.4] {$\ c$} (m-2-1)
    (m-2-3) edge[bend right=14] node[above=-0.1,pos=0.4] {$c, \lambda$} (m-2-1)
    (m-2-1) edge node[auto,swap] {$a$} (m-2-2)
    (m-2-3) edge node[auto] {$b$} (m-2-2)
    (m-3-1) edge node[auto,swap] {$a$} (m-3-2)
    (m-3-3) edge node[auto] {$b$} (m-3-2)
    (m-3-3) edge[bend right=14] node[auto,swap] {$c$} (m-3-1)
    (m-2-4) edge (m-3-4)
    ;
  \end{tikzpicture}
 \endpgfgraphicnamed
\]
In other words, $M$ is given by the three matrices $A_a=\begin{psmallmatrix} 1 & 0 \\ 0 & 1 \end{psmallmatrix}$, $A_b=\begin{psmallmatrix}1 & 0 \\0 & 1\end{psmallmatrix}$ and $A_c= \begin{psmallmatrix}\lambda & 0\\1 & \lambda\end{psmallmatrix}$. We consider the quiver Grassmannian $\Gr_\ue(M_\lambda)$ for dimension vector $\ue=(1,2,1)$. First note that the ambient product Grassmannian is 
\[
 \prod\Gr(e_q,d_q) \ = \ \Gr(1,2)\times\Gr(2,2)\times\Gr(1,2) \ = \ \P^1\times\P^1
\]
with Pl\"ucker coordinates $[\Delta_1:\Delta_4|\Delta_3:\Delta_6]$. There are no classical Pl\"ucker relations and there is precisely one quiver Pl\"ucker relation, which is
\[
 E(c,\emptyset,\{1,4\}) \ = \ \lambda\, \Delta_4\, \Delta_3 \, - \,\lambda\, \Delta_1\, \Delta_6 \, + \, \Delta_1 \, \Delta_3 \ = \ 0.
\]
From this, we see that $\Gr_\ue(M_\lambda)$ forms a flat family over $\C$ with respect to the parameter $\lambda$. While its fibres over $\lambda\neq 0$ are smooth quadrics, which are isomorphic to $\P^1$, the fibre over $\lambda=0$ is the transversal intersection of two projective lines in a point. 

This example exhibits a family of quiver Grassmannians over $\Spec \C[\lambda]$ associated with representations $M_\lambda$ of constant dimension vector for which the Euler characteristic of the general fibre $\P^1$ is $2$, but which assumes Euler characteristic $3$ in its special fibre at $\lambda=0$.

\subsection*{Example 3 (Elliptic curve and two projective lines)}
Let $Q$ be the generalized Kronecker with arrows $a$, $b$, $c$ and $d$. Let $M$ be the representation of dimension vector $(3,4)$ that is given by the following matrices:
\[
 A_a \ = \ \begin{psmallmatrix}
            1 & 0 & 0 \\
            0 & 0 & 0 \\
            0 & 0 & 0 \\
            0 & 0 & 1 \\
           \end{psmallmatrix}, 
     \qquad
 A_b \ = \ \begin{psmallmatrix}
            0 & 1 & 0 \\
            0 & 0 & 0 \\
            0 & 0 & 1 \\
            0 & 0 & 0 \\
           \end{psmallmatrix},
     \qquad
 A_c \ = \ \begin{psmallmatrix}
            0 & 0 & 0 \\
            0 & 0 & 1 \\
            0 & 0 & 0 \\
            1 & 0 & 0 \\
           \end{psmallmatrix}, 
     \qquad
 A_d \ = \ \begin{psmallmatrix}
            1 & 0 & 0 \\
            1 & 0 & 0 \\
            0 & 1 & 0 \\
            0 & 0 & 0 \\
           \end{psmallmatrix}.
\]
The associated coefficient quiver $\Gamma$ is as follows where we draw the different arrows separately for better readability:
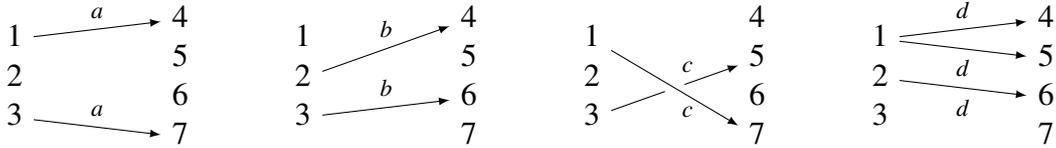
\[
 \beginpgfgraphicnamed{tikz/fig3}
  \begin{tikzpicture}[>=latex]
   \matrix (m) [matrix of math nodes, row sep=-0.5em, column sep=4em, text height=1ex, text depth=0ex]
    {        &    4    \\   
          1  &         \\   
             &    5    \\   
          2  &         \\   
             &    6    \\   
          3  &         \\   
             &    7    \\   
    };
   \path[->,font=\scriptsize]
    (m-2-1) edge node[above] {$a$} (m-1-2)
    (m-6-1) edge node[above] {$a$} (m-7-2)
    ;
  \end{tikzpicture}
 \endpgfgraphicnamed
 \qquad
 \beginpgfgraphicnamed{tikz/fig4}
  \begin{tikzpicture}[>=latex]
   \matrix (m) [matrix of math nodes, row sep=-0.5em, column sep=4em, text height=1ex, text depth=0ex]
    {        &    4    \\   
          1  &         \\   
             &    5    \\   
          2  &         \\   
             &    6    \\   
          3  &         \\   
             &    7    \\   
    };
   \path[->,font=\scriptsize]
    (m-4-1) edge node[above] {$b$} (m-1-2)
    (m-6-1) edge node[above] {$b$} (m-5-2)
    ;
  \end{tikzpicture}
 \endpgfgraphicnamed
 \qquad
 \beginpgfgraphicnamed{tikz/fig5}
  \begin{tikzpicture}[>=latex,back line/.style={},cross line/.style={preaction={draw=white, -,line width=6pt}}]
   \matrix (m) [matrix of math nodes, row sep=-0.5em, column sep=4em, text height=1ex, text depth=0ex]
    {        &    4    \\   
          1  &         \\   
             &    5    \\   
          2  &         \\   
             &    6    \\   
          3  &         \\   
             &    7    \\   
    };
   \path[->,font=\scriptsize]
    (m-6-1) edge [back line] node[above,pos=0.6] {$c$} (m-3-2)
    (m-2-1) edge [cross line] node[below,pos=0.6] {$c$} (m-7-2)
    ;
  \end{tikzpicture}
 \endpgfgraphicnamed
 \qquad
 \beginpgfgraphicnamed{tikz/fig6}
  \begin{tikzpicture}[>=latex]
   \matrix (m) [matrix of math nodes, row sep=-0.5em, column sep=4em, text height=1ex, text depth=0ex]
    {        &    4    \\   
          1  &         \\   
             &    5    \\   
          2  &         \\   
             &    6    \\   
          3  &         \\   
             &    7    \\   
    };
   \path[->,font=\scriptsize]
    (m-2-1) edge node[above] {$d$} (m-1-2)
    (m-2-1) edge node[below] {$d$} (m-3-2)
    (m-4-1) edge node[below] {$d$} (m-5-2)
    ;
  \end{tikzpicture}
 \endpgfgraphicnamed
\]
We consider the quiver Grassmannian $\Gr_\ue(M)$ for dimension vector $\ue=(1,3)$, which embeds into the product Grassmannian
\[
 \Gr(1,3)\times\Gr(3,4) \ = \ \Bigl\{ [ \, \Delta_1 : \Delta_2 :\Delta_3 \, \bigr| \, \Delta_{456} : \Delta_{457} : \Delta_{467} : \Delta_{567} \, \bigr] \Bigr\} \ \simeq \ \P^2\times\P^3.
\]
There are no classical Pl\"ucker relations for $\Gr(1,3)\times\Gr(3,4)$. The quiver Pl\"ucker relations are as follows:
\[
\begin{array}{lllll}
 E_I(a,\emptyset,\{4,5,6,7\}) & = & \Delta_3 \, \Delta_{456} \ - \ \Delta_1 \, \Delta_{567}                                  & = & 0, \\
 E_I(b,\emptyset,\{4,5,6,7\}) & = & \Delta_3 \, \Delta_{457} \ - \ \Delta_2 \, \Delta_{567}                                  & = & 0, \\
 E_I(c,\emptyset,\{4,5,6,7\}) & = & \Delta_3 \, \Delta_{467} \ - \ \Delta_1 \, \Delta_{456}                                  & = & 0, \\
 E_I(d,\emptyset,\{4,5,6,7\}) & = & \Delta_1 \, \Delta_{567} \ - \ \Delta_1 \, \Delta_{467} \ + \ \Delta_2 \, \Delta_{457}   & = & 0. \\
\end{array}
\]
On the locus where $\Delta_3$ and $\Delta_{567}$ are invertible, we can rewrite the first three equations as
\[
 \Delta_{456} \ = \ \Delta_{567} \, \Delta_1 \, \Delta_3^{-1}, \qquad \Delta_{457} \ = \ \Delta_{567} \, \Delta_2 \, \Delta_3^{-1}, \qquad \Delta_{467} \ = \ \Delta_{567} \, \Delta_1 \, \Delta_3^{-1},
\]
and replace $\Delta_{456}$, $\Delta_{457}$ and $\Delta_{467}$ in the latter equation, which becomes, after multiplication with $\Delta_3^2\Delta_{567}^{-1}$,
\[
 \Delta_2^2 \, \Delta_3 \ - \ \Delta_1^3 \ + \ \Delta_1 \, \Delta_3^2 \ = \ 0.
\]
This equation describes the projective closure of the locus where $\Delta_3$ and $\Delta_{567}$ are invertible, which is an elliptic curve.

It is easily verified that $\Delta_3$ and $\Delta_{567}$ vanish simultaneously, and we find the following additional irreducible components of $\Gr_\ue(M)$ in the locus $\Delta_3=\Delta_{567}=0$. The former two equations become trivial, and we are left with
\[
 \Delta_1 \, \Delta_{456} \ = \ 0 \quad \text{and} \quad \Delta_2 \, \Delta_{457} \ - \ \Delta_1 \, \Delta_{467} \ = \ 0.
\]
These equations describe the union of the two smooth rational lines 
\[
 \Bigr\{ \bigr[ \, 0:1:0 \, \bigr| \, \Delta_{456}:0:\Delta_{467}:0 \, \bigr] \Bigr\} \quad \text{and} \quad \Bigr\{ \bigr[ \, \Delta_1:\Delta_2:0 \, \bigr| \, 0:\Delta_{457}:\Delta_{467}:0 \, \bigr] \, \Bigr| \, \Delta_2 \Delta_{457} = \Delta_1 \Delta_{467} \, \Bigr\}.
\]
We conclude that $\Gr_\ue(M)$ is the union of an elliptic curve with two projective lines.


\begin{small}

\bibliographystyle{plain}

\begin{thebibliography}{1}

\bibitem{cr}
Philippe Caldero and Markus Reineke.
\newblock On the quiver {G}rassmannian in the acyclic case.
\newblock {\em J. Pure Appl. Algebra}, 212(11):2369--2380, 2008.

\bibitem{Hille15}
Lutz Hille.
\newblock Moduli of representations, quiver {G}rassmannians, and {H}ilbert schemes.
\newblock Preprint, \arxiv{1505.06008}, 2015.

\bibitem{L15}
Oliver Lorscheid.
\newblock Schubert decompositions for quiver {G}rassmannians of tree modules.
\newblock {\em Algebra Number Theory}, 9(6):1337--1362, 2015.
\newblock With an appendix by Thorsten Weist.

\bibitem{LW15a}
Oliver Lorscheid and Thorsten Weist.
\newblock Quiver {G}rassmannians of type {$\widetilde D_n$}. {P}art 1: {S}chubert systems and decompositions into affine spaces.
\newblock Accepted by {\em Memoirs of the AMS}, \arxiv{1507.00392}, 2015.

\bibitem{LW17}
Oliver Lorscheid and Thorsten Weist.
\newblock Representation type by Euler characteristics and singularities of quiver Grassmannians.
\newblock Preprint, \arxiv{1706.00860}, 2017.

\bibitem{Reineke12}
Markus Reineke.
\newblock Every projective variety is a quiver {G}rassmannian.
\newblock {\em Algebr. Represent. Theory}, 16(5):1313--1314, 2013.

\end{thebibliography}

\end{small}

\end{document}